\documentclass[a4paper,reqno,oneside,11pt]{amsart}

\usepackage[T1]{fontenc}
\usepackage{slantsc}

\usepackage{amsthm,amsmath,amssymb}
\usepackage{enumitem}
\usepackage{framed,graphics,color}

\usepackage[top=2cm, bottom=2cm, outer=3cm, inner=3cm,
            heightrounded, marginparwidth=3.5cm, marginparsep=1cm]{geometry}

\usepackage{setspace}
\usepackage[numbers,sort&compress]{natbib}

\usepackage{hyperref}
\usepackage{xspace}

\newtheorem{theorem}{Theorem}[section]

\newtheorem{lemma}[theorem]{Lemma}
\newtheorem{claim}[theorem]{Claim}
\newtheorem{corollary}[theorem]{Corollary}

\newtheorem{conjecture}[theorem]{Conjecture}

\newif\ifnotescoauthor
\notescoauthortrue 
\newcommand{\comment}[2]{\ifnotescoauthor \begin{framed} $\blacktriangleright$
{\sf #1} $\blacktriangleleft$ \ifx&#2&\else\medskip\hrule\medskip
#2\fi\end{framed}\fi}

\newcommand{\Prob}{{\mathbb{P}}}
\newcommand{\Bin}{{\textup{Bin}}}

\newcommand{\eps}{\varepsilon}

\newcommand{\rt}{\right}
\newcommand{\lt}{\left}
\newcommand{\cH}{\mathcal{H}}
\newcommand{\cS}{\mathcal{S}}
\newcommand{\cP}{\mathcal{P}}
\newcommand{\cHn}{\mathcal{H}^{\mathrm{N}}}
\newcommand{\SB}[1]{$(#1)$-SpookyBalance}

\newcommand{\EMAIL}[1]{  \textit{E-mail}: \texttt{#1}} 
\def\subset{\subseteq}

\def\itm#1{\rm ({#1})} 
\def\itmit#1{\itm{\it #1\,}} 
\def\rom{\itmit{\roman{*}}} 
\def\abc{\itmit{\alph{*}}}

\newcommand{\oldqed}{}
\def\endofClaim{\hfill\scalebox{.6}{$\Box$}}
\newenvironment{claimproof}[1][Proof]{
  \renewcommand{\oldqed}{\qedsymbol}
  \renewcommand{\qedsymbol}{\endofClaim}
  \begin{proof}[#1]
}{
  \end{proof}
  \renewcommand{\qedsymbol}{\oldqed}
}

\usepackage[mathlines]{lineno}
\usepackage{etoolbox} 

\newcommand*\linenomathpatch[1]{%
   \expandafter\pretocmd\csname #1\endcsname {\linenomath}{}{}%
   \expandafter\pretocmd\csname #1*\endcsname{\linenomath}{}{}%
   \expandafter\apptocmd\csname end#1\endcsname {\endlinenomath}{}{}%
   \expandafter\apptocmd\csname end#1*\endcsname{\endlinenomath}{}{}%
 }
\newcommand*\linenomathpatchAMS[1]{%
    \expandafter\pretocmd\csname #1\endcsname {\linenomathAMS}{}{}%
    \expandafter\pretocmd\csname #1*\endcsname{\linenomathAMS}{}{}%
    \expandafter\apptocmd\csname end#1\endcsname {\endlinenomath}{}{}%
    \expandafter\apptocmd\csname end#1*\endcsname{\endlinenomath}{}{}%
}

\expandafter\ifx\linenomath\linenomathWithnumbers
\let\linenomathAMS\linenomathWithnumbers
\patchcmd\linenomathAMS{\advance\postdisplaypenalty\linenopenalty}{}{}{}
\else
\let\linenomathAMS\linenomathNonumbers
\fi

\linenomathpatchAMS{gather}
\linenomathpatchAMS{multline}
\linenomathpatchAMS{align}
\linenomathpatchAMS{alignat}
\linenomathpatchAMS{flalign}
\linenomathpatch{equation}

\title{Making spanning graphs}

  \author[P. Allen]{Peter Allen*}
  \author[J. B\"ottcher]{Julia B\"ottcher*}

  \thanks{%
    * %
    Department of Mathematics, London School of Economics, Houghton Street,
    London WC2A 2AE, U.K.
    \EMAIL{\{p.d.allen|j.boettcher\}@lse.ac.uk}%
  }

  \author[Y. Kohayakawa]{Yoshiharu Kohayakawa\dag}

  \thanks{%
    \dag\
    Instituto de Matem\'atica, Estat\'{\i}stica e Ci\^encia da
    Computa\c c\~ao, Universidade de S\~ao Paulo, Rua do Mat\~ao 1010,
    05508--090~S\~ao Paulo, Brazil.  \EMAIL{yoshi@ime.usp.br}%
 }
  \author[H. Naves]{Humberto Silva Naves\ddag}
  \thanks{%
    \ddag\
    Institute for Mathematics and its Applications, University of
    Minnesota, Minneapolis, MN 55455, USA.
    \EMAIL{hnaves@ima.umn.edu}%
  }
  \author[Y. Person]{Yury Person\S}
  \thanks{%
    \S\
    TU Ilmenau, Weimarer Str.~25, 98684 Ilmenau, Germany.    
    \EMAIL{yury.person@tu-ilmenau.de}%
  }

  \thanks{YK is partially supported by
    FAPESP (2023/03167-5) and CNPq (407970/2023-1, 420838/2025-2,
    315258/2023-3).  YP is supported by the Carl Zeiss
    Foundation and by DFG (project no. 447645533).  The cooperation of the
    authors was supported by a joint CAPES-DAAD project
    (415/ppp-probral/po/D08/11629, Proj.~no.~333/09; 57350402).  This
    research was supported in part by the Institute for Mathematics
    and its Applications with funds provided by the National Science
    Foundation.  This research was partially supported by CAPES
    (Finance Code 001).
  }

\date{}

\begin{document}

\begin{abstract}
 We prove that for each $D\ge 2$ there exists $c>0$ such that whenever
 $b\le c\big(\tfrac{n}{\log n}\big)^{1/D}$, in the $(1:b)$ Maker-Breaker
 game played on $E(K_n)$, Maker has a strategy to guarantee claiming a
 graph~$G$ containing copies of all graphs $H$ with $v(H)\le n$ and
 $\Delta(H)\le D$. We show further that the graph~$G$ guaranteed by this strategy also contains
copies of any graph $H$ with bounded maximum degree and degeneracy at most
$\tfrac{D-1}{2}$. 
This lower bound on the threshold bias is sharp up to the $\log$-factor when $H$ consists of $\tfrac{n}{3}$ vertex-disjoint triangles or $\tfrac{n}{4}$ vertex-disjoint $K_4$-copies.
\end{abstract}
\maketitle

\section{Introduction}
\label{sec:introduction}

We study Maker-Breaker games on graphs, where the winning sets consist
of spanning subgraphs. So far this has been studied only in very
special cases such as spanning trees, Hamilton cycles and
clique-factors.  Here we consider general bounded-degree subgraphs.

Given integers $m$ and $b$, the $(m:b)$ \emph{Maker-Breaker game} played on a finite ground set~$X$ with \emph{winning sets} $\cS\subset\cP(X)$ is played as follows. In each round, first Maker claims up to $m$ unclaimed elements of~$X$. Then Breaker claims up to $b$ unclaimed elements of~$X$. If at some time Maker claims all elements of some $S\in\cS$, then Maker wins the game; otherwise Breaker wins. Since this is a finite draw-free game, either Maker or Breaker has a strategy which is guaranteed to win (see for example~\cite{GaleStewart}).
We remark that sometimes the game is defined with Breaker playing first
rather than Maker, but the results of this paper are insensitive to such a change.

We concentrate on the case $m=1$, and refer to $b$ as the \emph{bias}
of the game. Since the ability to claim more elements per round can
only help Breaker, for any game there is trivially a \emph{threshold
  bias}, defined to be the largest $b$ such that Maker has a winning
strategy. Note that for some games we have $b=0$, for example any game
in which the winning sets consist of pairwise disjoint sets of size
two. We are interested in games played on the ground set $E(K_n)$ in
which Maker's aim is to create, or ``claim'', a copy of a given graph,
where by a \emph{copy} of a graph~$H$ we mean a not necessarily
induced subgraph isomorphic to~$H$.  (Our results will in fact show
that Maker can simultaneously claim a copy of \textit{every} graph in
certain families of graphs.)

Games of this type received considerable attention. In particular the
\emph{connectivity game}, in which Maker's aim is to create any spanning
tree in $K_n$, was studied in~\cite{beck85:_random,ChvErd,GebSza}. More recently,
the Hamiltonicity game, in which Maker's aim is to create a Hamiltonian
graph, was investigated in~\cite{beck85:_random,BP82,ChvErd,Kriv11,KS08}. In both
cases the threshold bias is known quite accurately, and in fact is roughly
the same as if Maker simply chose $\tfrac{1}{b+1}\binom{n}{2}$ random
edges. Since Maker and Breaker do not in general play randomly, this is
rather surprising, and becomes more so as one observes that a similar
\emph{random graph intuition}
%
%
holds for many other graph games. 

The first result really harnessing this intuition is due to Bednarska and \L uczak in
\cite{BedLuc}. Given a graph $H$ on at least three vertices we define
\[
  m(H)=\max \lt\{\frac{|E(F)|-1}{|V(F)|-2}: F\subseteq H \text{ and }
      |V(F)|\ge 3\rt\}.
\]
Bednarska and \L uczak proved that the critical bias for the game
whose target sets are copies of $H$ is of order $\Theta\lt(n^{1/m(H)}\rt)$.


\begin{theorem}[Theorem 2 in \cite{BedLuc}]
  \label{thm:BedLuc}
  For every graph $H$ which contains a cycle there exists a constant
  $c_0$ such that for every sufficiently large integer $n$ and $b\le
  c_0n^{1/m(H)}$ the following holds. Consider the $(1:b)$ Maker-Breaker game
  played on the edges of $K_n$ with the copies of $H$ in $K_n$ as winning sets.
  Maker has a random strategy for this game that succeeds with
  probability $1-o(1)$ against any strategy of Breaker.
\end{theorem}

Maker uses the following simple random strategy.  She
\textit{secretly} picks a random ordering of~$E(K_n)$, and in her
$r$-th round claims the $r$-th edge in her ordering.  In case that
edge has already been claimed, she claims an arbitrary unclaimed edge
instead.  Notice that the first $m$ edges in Maker's ordering
correspond to the random graph $G(n,m)$, the random graph chosen
uniformly at random among all graphs with $m$ edges. It is well known
that $G(n,m)$ behaves in many respects as $G(n,p)$ with
$p=m/\binom{n}{2}$.  One then checks that, with this strategy, Maker
with high probability successfully claims most of a random graph
$G(n,p)$, where $p$ is a very small constant multiple of $b^{-1}$.
Bednarska and {\L}uczak proved a \emph{global resilience} result which
states that, if $p\gg n^{-1/m(H)}$, then the random graph
$\Gamma=G(n,p)$ typically has the property that any $G\subseteq\Gamma$
which contains most of the edges of~$\Gamma$ also contains a copy of
$H$, justifying that Maker wins the $H$-game with bias~$b$ with high
probability.

Since a Maker-Breaker game is a deterministic game, it follows that if
Maker has a random strategy that works with non-zero probability
against any given strategy of Breaker, then the game is Maker's win
(otherwise Maker's strategy should work with probability zero against
Breaker's winning strategy; arguments of this nature go back at least
to~\cite{spencer94:_random}). Thus there exists a deterministic
strategy which wins the $H$-game for Maker at the same bias, though
the proof does not give any hint of what this strategy might be.

It is easy to see that for some games on $E(K_n)$
the Bednarska--{\L}uczak strategy will not work close to the threshold
bias. For example, it is not a good strategy in the connectivity game:
Breaker can choose all his edges at a fixed vertex, in which case Maker
must claim an edge at that vertex in the first $\tfrac{n-1}{b}$ turns or it
will be isolated and she will lose. This is an unlikely event whenever
$b=\omega(1)$, yet the threshold bias for the connectivity game is
$\Theta\big(\tfrac{n}{\log n}\big)$ (for the best bound see Gebauer and
Szab\'o~\cite{GebSza}). 

Overcoming this particular obstacle, Ferber, Krivelevich and
Naves~\cite{FKN} gave a different randomised strategy which, provided
the bias is, say, $o\big(\tfrac{n}{\log n}\big)$, ensures that Maker
claims most edges of an \textit{underlying random graph} at each
vertex.  We remark that games in which the winning criterion has to do
with degrees had already been considered; see, e.g., \cite{balogh09}
for sharp results on the so called \textit{biased degree game}.
However, here, we are also concerned with a certain underlying random
graph~$G(n,p)$, which Maker generates while playing the game.  At the
end of the game, Maker has to present a graph~$\Gamma$ distributed
as~$G(n,p)$, and she must have claimed a very large proportion of the
edges of~$\Gamma$ at every vertex.  A little more formally, for any
fixed Breaker strategy, Maker secretly draws a random
graph~$\Gamma=G(n,p)$ with $\frac{\log n}{n}\ll p\ll b^{-1}$, which,
with high probability, is such that Maker's strategy succeeds in
claiming a subgraph $G$ of $\Gamma$ with minimum degree nearly~$pn$.
The advantage of this strategy is that now a \emph{local resilience}
result on $G(n,p)$ gives Maker a winning strategy in the
``corresponding'' $(1:b)$ game, and several such local resilience
results were known at the time~\cite{FKN} was written.  This reproved,
for example, that the threshold bias in the Hamiltonicity game is
$\Theta\big(\tfrac{n}{\log n}\big)$, as first established by
Beck~\cite{beck85:_random}.  As is well known, this result has been
refined to the asymptotically optimal value of $(1-o(1))n/\log n$ by
Krivelevich~\cite{Kriv11}.

The randomised strategy of~\cite{FKN} is, roughly, to reveal edges of
$G(n,p)$ not arbitrarily but in an order responsive to Breaker's
actions, revealing edges preferentially at vertices where Breaker has claimed more edges.  However, as it turns out, this strategy does not allow Maker
to win, for example, the $K_3$-factor game with an optimal bias
because Breaker can easily focus on a fixed vertex~$v$ and claim all
edges that appear in the neighbourhood of~$v$.  That way~$v$ will not
be in any triangle and hence there will not be any $K_3$-factor.\footnote{When Maker claims an edge $vu$ at $v$, Breaker will in his next $\le pnb^{-1}$ turns try to claim all the edges $ux$ such that Maker previously claimed $vx$. Heuristically, we expect Breaker completes these turns before Maker claims another edge at $v$. Maker will claim roughly $pn$ edges at $v$ in total, and so she will have roughly $p^2n^2b^{-1}$ turns in which she has a chance to claim an edge making a triangle with $v$. Heuristically, a uniform choice of edge in each turn would pick such an edge with probability at most $pn\binom{n}{2}^{-1}\approx \tfrac{2p}{n}$, and we claim the strategy of~\cite{FKN} will not (significantly) prioritise such edges. Thus the probability Maker will get a triangle at $v$ is bounded by $2p^2n^2 b^{-1}\cdot\tfrac{2p}{n}= 4p^3n b^{-1}$. Since $p\ll b^{-1}$ we see Breaker likely wins already when $b\gg n^{1/4}$.}

Our
main result overcomes this obstacle.

\subsection*{Our results}

Our main result states that if $p\gg n^{-1/2}$ and $b\ll p^{-1}$, in
the $(1:b)$ Maker-Breaker game played on $K_n$, Maker has a randomised
strategy which, on the one hand generates a random
graph~$\Gamma=G(n,p)$, and, on the other hand, guarantees that she
will claim close to~$pn$ edges of~$\Gamma$ at each
vertex~$v\in V(\Gamma)$, and will also claim nearly $p^3n^2/2$ edges
of~$\Gamma$ in each $\Gamma$-neighbourhood of~$v$ for every
vertex~$v\in V(\Gamma)$.

\begin{theorem}\label{thm:makeGnp} 
  For every~$n$, $\eps\in(0,1)$, $p\ge10^8\eps^{-2}n^{-1/2}$, and
  $b\le 10^{-24}\eps^6 p^{-1}$ the following holds.  There is a
  randomised strategy for Maker in the $(1:b)$ Maker-Breaker game
  played on~$K_n$ that produces a random graph $\Gamma=G(n,p)$ and is
  such that, with high probability, lets Maker claim a subgraph~$G$
  of~$\Gamma$ with $\delta(G)\ge(1-\eps)pn$ and
  $e\big(G[N_\Gamma(v)]\big)\ge (1-\eps)p^3n^2/2$.
\end{theorem}

It is easy to see that, up to constant factors, this result is sharp
in the following sense.  Whatever strategy Maker plays, if
$b\gg p^{-1}$ then Breaker can simply choose edges at a fixed vertex
until none remain, and so Maker claims at most $\tfrac{n}{b+1}\ll pn$
edges at that vertex.  Similarly, if $p\ll n^{-1/2}$ and
$b=\Omega(p^{-1})$, then $b\gg n^{1/2}$, and, as proved in the seminal
paper by Chv\'atal and Erd\H os~\cite{ChvErd}, Breaker has a strategy
that stops Maker from creating any triangle.  Of course, this means
that, for such~$p$ and~$b$, Breaker can force every vertex~$v$ to be
such that $e\big(G[N_G(v)]\big)=0$, so $e\big(G[N_\Gamma(v)]\big)$ will also be small.

Following the spirit of~\cite{frieze05:_jumbl}, in which the authors
obtain several results on positional games by showing how Maker can
achieve pseudo-random graphs, we deduce from Theorem~\ref{thm:makeGnp}
upper bounds on the critical bias of Maker-Breaker games whose targets
are spanning subgraphs of~$K_n$ of low maximum degree and low
degeneracy.  We do this by making use of results from~\cite{SparseBU}
to show that if $\Gamma$ is a typical sample from $G(n,p)$ then any
subgraph $G$ of $\Gamma$ satisfying the conditions in
Theorem~\ref{thm:makeGnp} contains all low-maximum-degree and
low-degeneracy graphs.

Recall that a graph is $D$-degenerate if all of its subgraphs contain
a vertex of degree at most~$D$.  Recall also that the existence of a
randomised strategy for Maker which wins with positive probability
implies that Breaker has no winning strategy.

\begin{corollary}\label{cor:makeuniv}
 For every~$D$ and~$\Delta$, there exists $c>0$ such that the following
 holds for any $b\le c\big(\tfrac{n}{\log n}\big)^{1/D}$. 
 In the $(1:b)$ Maker-Breaker game played on $E(K_n)$, Maker has a
 deterministic strategy to claim a graph~$G$ which contains copies of all
 graphs $H$ on at most $n$ vertices such that
 \begin{enumerate}[label=\rom]
   \item $\Delta(H)\le D$, or
   \item $\Delta(H)\le\Delta$ and $H$ is $\tfrac{D-1}{2}$-degenerate.
 \end{enumerate}
\end{corollary}

In particular, there is $c>0$ such that Maker can obtain a
$K_3$-factor provided Breaker's bias is at most
$c\big(\tfrac{n}{\log n}\big)^{1/2}$, and a $K_4$-factor if the bias
is at most $c\big(\tfrac{n}{\log n}\big)^{1/3}$.  These bounds are
optimal up to the $\log$-factor.  Indeed, as proved by Chv\'atal and
Erd\H os~\cite{ChvErd}, if $b\geq2n^{1/2}$, then Breaker can prevent
Maker from claiming even a single triangle (for improvements of the
constant~$2$ in the lower bound of~$b$ above, see~\cite{glazik22,
  balogh11:_chvat_erdos}).  The threshold bias for the $K_r$-factor
game, for~$r\geq4$, has been established by Liebenau and
Nenadov~\cite{liebenau22}: it is~$n^{2/(r+2)}$ up to multiplicative
constants. 

We do not believe Corollary~\ref{cor:makeuniv} is optimal up to the $\log$-factor in general,
because we expect that the results of~\cite{SparseBU} can be improved.
We suspect that it is possible (though perhaps technically hard) to
remove the $\log$-factors from those results, but we also believe that
the polynomial term can in general be better.  In the spirit of
inciting further research, however, we conjecture that the strategy of
Theorem~\ref{thm:makeGnp} is optimal up to a constant for these
problems.

\begin{conjecture}
  For each $\Delta$, there is a constant $C$ for which the following
  holds.  Let~$H$ be an $n$-vertex graph with $\Delta(H)\le\Delta$,
  and suppose that Maker has a winning strategy to claim a copy of~$H$
  in the $(1:Cb)$ game on~$E(K_n)$.  Then, in the
$(1:b)$ game on $E(K_n)$, playing Maker's strategy presented in the
proof of Theorem~\ref{thm:makeGnp}, Maker is able to claim a copy
of~$H$ with high probability.
\end{conjecture}

We remark that if Maker plays the strategy of
Theorem~\ref{thm:makeGnp} using the random graph $G(n,p)$ when
$\tfrac{\log n}{n}\ll p\ll n^{-1/2}$, while she will not necessarily
claim edges in the neighbourhood of all vertices, she will at least
succeed in making a subgraph of $G(n,p)$ with minimum degree close to
$pn$. Similarly, if $n^{-2}\ll p\ll\tfrac{\log n}{n}$,
she will not achieve this minimum degree, but she will at least claim
most of the edges of~$G(n,p)$.

Going back to the problem of finding an optimal result for the
Maker-Breaker game for spanning target graphs, we mention that, for
$\Delta\geq3$, Liebenau and Nenadov conjecture that
$K_{\Delta+1}$-factors are essentially the hardest target graphs among
all graphs of maximum degree~$\Delta$
(see~\cite[Conjecture~1.3]{liebenau22}).

Finally, we observe that while we prove the existence of a deterministic winning strategy for Maker in the game of Corollary~\ref{cor:makeuniv}, our proof gives no hint of what that strategy might be. It would be interesting to find such deterministic winning strategies, even for a simple case such as the $K_3$-factor game.

\subsection*{Organisation}

In Section~\ref{sec:auxiliary} we collect some tools that we need in our
proofs. In Section~\ref{sec:outline} we outline the proof of
Theorem~\ref{thm:makeGnp}, for which we need an
auxiliary game that is defined and analysed 
in Section~\ref{sec:boxgame}. The details of the proof
of Theorem~\ref{thm:makeGnp} are then provided in Section~\ref{sec:maker}. 
In Section~\ref{sec:resilience} we prove a resilience result, which together
with Theorem~\ref{thm:makeGnp} implies
Corollary~\ref{cor:makeuniv}. 

\subsection*{Notation}

We sometimes write $a\pm b$ for a quantity~$x$ for which $a-b\leq
x\leq a+b$ holds.  We omit floor and ceiling symbols when they are not
essential.

\section{Preliminaries}
\label{sec:auxiliary}

In this section we present some auxiliary results that will be used
throughout the paper. We use extensively the following well-known
bound on the tails of the binomial distribution
due to Chernoff (see, for instance, \cite{JLRbook}).

\begin{lemma}
  \label{lem:chernoff}
  If $X \sim \mathrm{\Bin}(n,p)$, then for every $0<a<\tfrac32$ we have
  \[\Prob\big[|X-np|>anp\big]<2\exp\Big(-\frac{a^2np}{3}\Big)\,.\]
\end{lemma}

We also use the following graph partitioning lemma, derived in~\cite[Lemma~2.4]{SparseBU} from the Hajnal-Szemer\'edi Theorem~\cite{HajSze70}.

\begin{lemma}\label{lem:nicepartition}
  Given any graph $F$ and a subset $Y$ of $V(F)$, there is a partition $V(F)=V_1\cup\dots\cup V_{8\Delta(F)}$ such that each part is
  independent, the parts differ in size by at most one, and the sets $V_i\cap Y$ differ in size by at most one.
\end{lemma}

Finally, we need to know that dense subgraphs of typical random graphs contain all sparse graphs. This is one of the main results of~\cite{SparseBU}. The statement we give is a weakening and combination of Lemmas~1.21 and~1.23 from that paper.
To state it, we need the following definition. We say that a pair $(X,Y)$ of disjoint vertex sets in a graph $G$ is \emph{$(\eps,d,p)$-regular} if for any $X'\subseteq X$ and $Y'\subseteq Y$, we have $e(X',Y')\ge(d-\eps)p|X'||Y'|$. 
For a graph~$\Gamma$, a vertex $v\in V(\Gamma)$, and a vertex set $A\subset
V(\Gamma)$ we also write $N_\Gamma(v;A)$ for the neighbours of~$v$ in~$A$, and we simply write $N_\Gamma(v)$ for $N_\Gamma(v;V(\Gamma))$.

\begin{lemma}\label{lem:simpBL}
  For all $r,\Delta,D\ge 2$ and $d>0$ there exist $\eps\in (0,1)$ and $C$ such that if
  $p\ge C(\log n/n)^{1/D}$, then the random graph $\Gamma=G(n,p)$ asymptotically almost
  surely has the following property.
  
  Let $G\subseteq\Gamma$ and $H$ be graphs on $n$
  vertices. Suppose that $V(G)=V_1\cup\cdots\cup V_{r}$ is a 
  partition with parts differing in size by at most one, and that $H$ has an $r$-colouring with colour classes
  $X_1,\ldots,X_{r}$ where $|V_i|=|X_i|$ for each $i$. Suppose that
  \begin{enumerate}[label=\rom]
   \item $\Delta(H)\le D$, or
   \item $\Delta(H)\le\Delta$, and $H$ has degeneracy at most $\tfrac{D-1}{2}$, and at least $\tfrac{1}{4D^2}|X_i|$ vertices of each $X_i$ have degree at most $2D$.
  \end{enumerate}
  Suppose
  furthermore that the following conditions hold.
  \begin{enumerate}[label=\abc]
    \item\label{simpBL:reg} For each $i<j$, the pair $(V_i,V_j)$ is
    $(\eps,d,p)$-regular in $G$.
    \item\label{simpBL:degree} For each $i\neq j$ and each $v\in V_i$, we
    have $\big|N_G(v;V_j)\big|\ge (d-\eps)p|V_j|$.
    \item\label{simpBL:1reg} For each $i\neq j\neq k$ and each $v\in V_i$, the
    pair $\big(N_\Gamma(v;V_j),V_k\big)$ is $(\eps,d,p)$-regular in $G$.
    \item\label{simpBL:2reg} For all distinct indices $i,j,k$ and all $v\in
    V_i$, the pair $\big(N_\Gamma(v;V_j),N_\Gamma(v;V_k)\big)$ is $(\eps,d,p)$-regular in $G$.
  \end{enumerate}
  Then $H$ is a subgraph of $G$.
\end{lemma}

\section{Outline of the strategy}
\label{sec:outline}

We now briefly describe Maker's strategy in Theorem~\ref{thm:makeGnp}.
Maker will not generate $\Gamma=G(n,p)$ at once.  Rather, in each
round, she will pick an edge of $K_n$ not yet chosen by Breaker, and
query whether this edge is in~$\Gamma=G(n,p)$ by flipping a $p$-biased
coin.  If the edge turns out to be in~$\Gamma$, she will claim it and
will put it in her graph~$G$; if not, she will pick another edge, and
so on, until she finds an edge of~$\Gamma$ to add to~$G$.  For
convenience, we suppose that Maker reveals to Breaker not only which
edges she claimed, but also which edges she queried. Note that under
this strategy Breaker never has an incentive to claim any edge which
Maker has queried and found not to be in~$\Gamma$.  We therefore can
and will assume that Breaker never does claim such edges.  When
Breaker claims an edge, we will defer the decision as to whether or
not that edge belongs to~$\Gamma=G(n,p)$.  We stop the game when all
the edges of~$K_n$ are either queried by Maker or claimed by Breaker.
However, formally speaking, Maker has to produce ``all
of''~$\Gamma=G(n,p)$, and hence, at the end of the game, she flips a
coin for every edge claimed by Breaker to ``complete'' her
subgraph~$G\subset\Gamma$ to~$\Gamma=G(n,p)$.

The strategic part of Maker's play is thus to determine the order in
which she queries edges of~$K_n$ to check whether or not they are
in~$\Gamma$.  The idea is as follows. In alternating turns, Maker will
try to avoid Breaker claiming too many edges at a given vertex (the
\emph{degree game}), or in a given vertex neighbourhood (the
\emph{neighbourhood game}), respectively. In both cases, she does this
by (privately) playing an auxiliary game called \emph{balance game}. A balance
game is played on the vertex set (for us, corresponding to $E(K_n)$)
of a hypergraph~$\cH$. The idea of
using balance games goes back to the box game of Chv\'atal and Erd\H{o}s~\cite{ChvErd}.
In this paper, following Gebauer and
Szab\'o~\cite{GebSza} and Ferber, Krivelevich and Naves~\cite{FKN},
Maker aims to claim almost a `fair share' of each
hyperedge. Since Breaker might choose to concentrate on only one of the
degree and neighbourhood games, in effect the individual games are
$(1:2b)$ Maker-Breaker games, and thus a `fair share' is a
$\tfrac{1}{1+2b}$-fraction of each hyperedge.

As is well known, the \emph{potential function strategy} of
Beck~\cite{Beck} is a good strategy for Maker to win games similar to box games. This
strategy assigns to each vertex of $\cH$ a score, and Maker claims the
highest-scoring unclaimed vertex in each round. However there are two
complicating factors in the games we consider. First, Maker typically
cannot claim the highest-scoring unclaimed edge, because she can claim
only edges of $\Gamma$. Second, in the neighbourhood game the hyperedges
change over time. This is the case because during the game the
$\Gamma$-neighbourhoods of the vertices are decided in stages, and
only their portion which Maker revealed so far by claiming edges is
known.

In Section~\ref{sec:boxgame} we define our balance
game, which takes into account these complicating factors, which we model as
the actions of a third player called \emph{Ghost}. We prove (in Lemma~\ref{lem:SB}) that a
modification of Beck's strategy allows Maker to win this game. In this
proof it turns out that we do not need to know that the Ghost's actions
result from a combination of the randomness of $\Gamma$ and Maker's own
actions. Instead, it suffices to require hyperedges to change over time only by
adding unclaimed vertices of $\cH$ (which correspond to edges of $K_n$),
and that this change has to take place immediately before Maker's next turn.

Analysing these two balance games, we will show that typically Breaker claims
only a tiny fraction of the edges of $K_n$ at each vertex and in each
vertex neighbourhood, before we stop the game. To understand this, recall that we stop the game when all edges of $K_n$ have been either claimed by Breaker or queried by Maker, and that Breaker has no incentive to claim edges Maker queried, so we assume he does not. Thus we expect Maker will query about $p^{-1}\gg b$ edges for every $1$ edge that she manages to claim.
Since by our assumption these edges claimed by Breaker were not revealed
in~$\Gamma$ during the game, with high probability~$\Gamma$ also only
contains a tiny fraction of these edges at each vertex and in each vertex
neighbourhood, proving Theorem~\ref{thm:makeGnp}.

Formally, the Maker-Breaker game ends only when every edge of $K_n$ is claimed by either Maker or Breaker, so once Maker has revealed all of $\Gamma$, she must still continue to play, claiming arbitrary available edges, until the game ends. By construction these arbitrarily chosen edges are not in $E(\Gamma)$, so do not affect the conclusion of Theorem~\ref{thm:makeGnp}.



\section{\texorpdfstring{A game with spirit \\ {\it In which a Ghost appears in
    the game, but in the end nothing really changes}}{A game with spirit}}
\label{sec:boxgame}

In this section we define and analyse the auxiliary balance game used in
the proof of Theorem~\ref{thm:makeGnp}.  This is a perfect information
game played by three players, Maker, Breaker and Ghost, and we call it
the \emph{\SB{m,b,V,e,\ell,M} game}. It is played on a vertex set
$V$. Initially a multihypergraph $\cH_0$ on $V$ is given, which we
call the \emph{underlying multihypergraph}, with $e$ hyperedges, some of
which may be the empty set and all of which have at most $M$
vertices. We refer to the vertices in~$V$ taken by Maker or Breaker as
\emph{claimed}, whereas the vertices taken by the Ghost are
\emph{haunted}. At the beginning all vertices in~$V$ are unclaimed and
unhaunted.  Now, in each round $r\ge 1$, play goes as follows. First,
the Ghost creates $\cH_r$ from $\cH_{r-1}$ by, for each hyperedge
of~$\cH_{r-1}$, adding arbitrarily some (possibly empty) set of
vertices of $V$ which are neither claimed nor haunted, subject to the
constraint that no hyperedge of $\cH_r$ has more than $M$ vertices. Next,
Maker chooses an unclaimed and unhaunted vertex. The Ghost either
allows her to claim the vertex she chooses, or itself \emph{haunts}
the vertex, stopping her from claiming the vertex.  The Ghost may
choose in addition to haunt (arbitrarily) other vertices as well.  In
each round, Maker continues to choose vertices until she has
claimed~$m$ vertices successfully or all vertices are either claimed
or haunted. Finally, Breaker chooses $b$ vertices which are neither
claimed nor haunted, and claims all of them. The game stops when all
vertices of $V$ are either claimed or haunted.

The game is won by Maker if at every round $r$ and for each hyperedge $S\in
E(\cH_r)$, the following holds. Let~$c$ be the number of vertices of~$S$
that are claimed. Then at least $\tfrac{mc}{m+b}-\ell$ vertices of $S$ are claimed by Maker.

In the remainder of this section we shall find, given $m$, $b$, $V$,
$e$ and $M$, a lower bound on $\ell$ such that Maker has a winning strategy
for the \SB{m,b,V,e,\ell,M} game. This strategy is a \emph{potential function
  strategy} with parameters $0<\lambda<1$ and $\tau>0$, which we now
describe. If Maker has claimed $X\subset V$ and Breaker has claimed
$Y\subset V$, we define the potential of a set $A\subseteq V$ to be 
\[\Phi(A,X,Y):=(1-\lambda)^{|A\cap X|}(1+\tau)^{|A\cap Y|}\,.\]
We define the potential of a multihypergraph $\cH$ on $V$ to be
\[\Phi(\cH,X,Y):=\sum_{S\in E(\cH)}\Phi(S,X,Y)\,,\]
and the effect of $v\in V$ (with respect to $\cH$) to be
\[\overline{\Phi}(v,\cH,X,Y):=\Phi(\cH,X,Y)-\Phi(\cH,X\cup\{v\},Y)\,.\]
Maker's strategy is simple: when Maker has a choice, if $\cH$ is the
current multihypergraph, and $X$ and $Y$ are the sets currently
claimed by Maker and Breaker respectively, then she chooses a vertex
$v$ which maximises $\overline{\Phi}(v,\cH,X,Y)$ among all unclaimed
and unhaunted vertices. We call this strategy the
\emph{$(\lambda,\tau)$-potential strategy}.

We remark that, if the Ghost never takes any action (i.e.\ allows
Maker to claim all the vertices she chooses, and does not add vertices
to edges of $\cH_0$), this game is a Balancer-Unbalancer game (see
e.g.~\cite[Section~16]{Beck}), and the potential function strategy is
well known to be good for Maker playing Balancer (see, in particular,
\cite[p.~300]{Beck}). The following lemma states that the Ghost's
actions do not have any material effect on the game---as one would
expect of a Ghost---and the proof is essentially the observation that
the standard proof for the Balancer-Unbalancer game goes through
unaffected by the actions of the Ghost.

\begin{lemma}\label{lem:SB}
  Let $m,\,b,\,e,\,M\in\mathbb{N}$ with $M\ge 9(m+b)\log e$ be given,
  and let~$\cH_0$ be a multihypergraph with~$e$ hyperedges on a vertex
  set~$V$.  Suppose 
  \[\ell\ge\frac{5mb}{m+b}\cdot \sqrt{\frac{M\cdot  \log e}{m+b}}\,.\]
  Then there exist $0<\lambda<1$ and~$\tau>0$ such that the
  $(\lambda,\tau)$-potential strategy is a winning strategy for Maker
  in the \SB{m,b,V,e,\ell,M} game.
\end{lemma}
\begin{proof}
 Given $m$, $b$, $e$ and $M$ we set
\[
  \lambda = \frac{1}{m}\sqrt{\frac{(m+b)\log e}{M}}
  \leq\frac1{3m}\leq\frac13\,,
\]
 and let $\tau$ be the unique positive solution of the equation $(1+\tau)^b=1+m\lambda$.

 Fix arbitrary strategies of Breaker and Ghost, and let Maker play the
 $(\lambda,\tau)$-potential strategy. Let $\cH_0$ be given and let
 $X_0=Y_0=\emptyset$. For each $r\ge 1$, let $X_r$ be the set of vertices
 claimed by Maker after round $r$, let $Y_r$ be the set of vertices claimed
 by Breaker after round $r$, let $Z_r$ be the haunted vertices after round
 $r$, and recall that~$\cH_r$ is the multihypergraph created by the Ghost in round $r$. Observe that $\Phi(\cH_0,X_0,Y_0)=e$.
 \begin{claim}\label{cl:SB:dec} $\Phi(\cH_r,X_r,Y_r)$ is monotone decreasing in $r$.
 \end{claim}
 \begin{claimproof}
 Fix $r\ge 1$. We aim to show  $\Phi(\cH_{r-1},X_{r-1},Y_{r-1})\ge \Phi(\cH_r,X_{r},Y_{r})$.
 First observe that $\Phi(\cH_{r-1},X_{r-1},Y_{r-1})=\Phi(\cH_r,X_{r-1},Y_{r-1})$ since the Ghost only adds unclaimed vertices to hyperedges. Suppose that in round $r$, Maker claims vertices $u_1,\dots,u_m$ and Breaker then claims vertices $v_1,\dots,v_b$. Since the effect of Maker claiming a vertex $u_i$ on the potential of any given hyperedge $S\in E(\cH_{r})$ is either nothing (if $u_i\not\in S$) or to multiply it by a factor $1-\lambda$ (if $u_i\in S$), it follows that the effect of unclaimed vertices decreases as Maker claims the vertices $u_1,\dots,u_m$ in turn. Let
 \[\overline{\Phi}_r=\max_{v\in V\setminus (X_r\cup Y_{r-1}\cup Z_r)}\overline{\Phi}(v,\cH_{r},X_r,Y_{r-1})\,.\]
 Since Maker chooses an unclaimed and unhaunted vertex maximising the
 effect, we have
 \begin{equation}
   \label{eq:2}
  \Phi(\cH_{r},X_r,Y_{r-1})-\Phi(\cH_{r},X_{r-1},Y_{r-1})\le
   -m\cdot\overline{\Phi}(v,\cH_{r},X_r,Y_{r-1})
 \end{equation}
 for each $v\in V\setminus (X_r\cup Y_{r-1}\cup Z_r)$, and so, in particular,
 \[\Phi(\cH_{r},X_r,Y_{r-1})-\Phi(\cH_{r},X_{r-1},Y_{r-1})\le -m\cdot\overline{\Phi}_r\,.\]
 
 Now by definition, for any $\cH$ and $X,Y$, and any
 $v\not\in X\cup Y$, we
 have
 \[\Phi(\cH,X,Y\cup\{v\})-\Phi(\cH,X,Y)=\tfrac{\tau}{\lambda}\overline{\Phi}(v,\cH,X,Y)\,.\]
 Furthermore, for any $Y'\supseteq Y$ and any $A\subseteq V$, we have
 \[\Phi(A,X,Y')\le(1+\tau)^{|Y'|-|Y|}\Phi(A,X,Y)\,.\]
 We conclude that
 for each $1\le i\le b$ we have
 \begin{align*}
   &\Phi\big(\cH_{r},X_r,Y_{r-1}\cup\{v_1,\dots,v_i\}\big)-\Phi\big(\cH_{r},X_r,Y_{r-1}\cup\{v_1,\dots,v_{i-1}\}\big)\\
   &\qquad\qquad=\tfrac{\tau}{\lambda}\overline{\Phi}\big(v_i,\cH_r,X_r,Y_{r-1}\cup\{v_1,\dots,v_{i-1}\}\big)
     \le\tfrac{\tau}{\lambda}(1+\tau)^{i-1}\overline{\Phi}(v_i,\cH_r,X_r,Y_{r-1})\\
   &\qquad\qquad\le\tfrac{\tau}{\lambda}(1+\tau)^{i-1}\overline{\Phi}_r \,.
 \end{align*}
 Summing over $1\le i\le b$, we have a geometric series, and so
 \[\Phi(\cH_{r},X_r,Y_r)-\Phi(\cH_{r},X_r,Y_{r-1})\le\frac{\tau}{\lambda}\cdot\frac{(1+\tau)^b-1}{(1+\tau)-1}\overline{\Phi}_r=\frac{(1+\tau)^b-1}{\lambda}\overline{\Phi}_r\,.\]
 At last, in view of~\eqref{eq:2}, we conclude that
 \[\Phi(\cH_{r},X_r,Y_r)-\Phi(\cH_{r},X_{r-1},Y_{r-1})\le \frac{(1+\tau)^b-1}{\lambda}\overline{\Phi}_r-m\cdot\overline{\Phi}_r= 0\,,\]
 where the final equality uses our choice of $\tau$ such that $(1+\tau)^b=1+m\lambda$.
 \end{claimproof}
 
 Now suppose that Maker does not win the \SB{m,b,V,e,\ell,M} game. Then there is some stage $r$ and hyperedge $S\in E(\cH_r)$ such that $|X_r\cap S|<\tfrac{m}{m+b}\big|(X_r\cup Y_r)\cap S\big|-\ell$. In particular, we have
 \[\Phi(S,X_r,Y_r)=(1-\lambda)^{|X_r\cap S|}(1+\tau)^{|Y_r\cap S|}\ge \big(\tfrac{1+\tau}{1-\lambda}\big)^\ell\cdot \big((1-\lambda)^m(1+\tau)^b\big)^{|(X_r\cup Y_r)\cap S|/(m+b)}\,.\]
 By choice of $\tau$, we have $(1+\tau)^b=1+m \lambda$, and we have
 \[(1-\lambda)^m(1+m \lambda)=(1-\lambda)^{m-1}\big(1+(m-1)\lambda-m\lambda^2\big)< (1-\lambda)^{m-1}\big(1+(m-1)\lambda\big)\,.\]
Thus we have $(1-\lambda)^m(1+\tau)^b<1$. In particular, this gives
 \[\Phi(S,X_r,Y_r)=(1-\lambda)^{|X_r\cap S|}(1+\tau)^{|Y_r\cap S|}\ge \big(\tfrac{1+\tau}{1-\lambda}\big)^\ell\cdot \big((1-\lambda)^m(1+\tau)^b\big)^{|S|/(m+b)}\,.\]
 Finally, because $S\in E(\cH_r)$ we have $|S|\le M$. Plugging in our chosen $\tau$, we have
 \begin{equation}\label{eq:finpot}
  \log \Phi(S,X_r,Y_r)\ge \ell\big(\tfrac{1}{b}\log(1+m\lambda)-\log(1-\lambda)\big)+\tfrac{M}{m+b}\big(m\log(1-\lambda)+\log(1+m\lambda)\big)
  \,.
 \end{equation}
 
Since $\log(1+m\lambda)\ge m\lambda-m^2\lambda^2/2$ and $\log(1-\lambda)<-\lambda$ hold, we thus have 
\[\tfrac{1}{b}\log(1+m\lambda)-\log(1-\lambda)>\tfrac{m}{b}\lambda-\tfrac{m^2}{2b}\lambda^2+\lambda=\tfrac{m+b}{b}\lambda-\tfrac{m^2\lambda^2}{2b}\ge\tfrac{m+b}{2b}\lambda\,,\]
where the final inequality uses $m\lambda\le \tfrac13$. 
 
 Similarly, we use $\log(1-\lambda)\ge -\lambda-\lambda^2/2-\lambda^3$ for $\lambda\in(0,1/3)$ to show
 \[m\log(1-\lambda)+\log(1+m\lambda)\ge -m\lambda-m\tfrac{\lambda^2}{2}-m\lambda^3+m\lambda-\tfrac{m^2\lambda^2}{2}>-\tfrac{3m^2\lambda^2}{2}\,. \]
 
 Plugging these approximations into~\eqref{eq:finpot}, and using the equalities $\lambda\ell=\tfrac{5b}{m+b}\log e$ and $m^2\lambda^2=\tfrac{m+b}{M}\log e$, we have
 \[\log \Phi(S,X_r,Y_r)> \ell\tfrac{m+b}{2b}\lambda-\tfrac{3m^2\lambda^2}{2}\cdot\tfrac{M}{m+b}\ge\tfrac{5}{2}\log e-\tfrac{3}{2}\log e=\log e\,.\]
 By Claim~\ref{cl:SB:dec} we have $e\ge\Phi(\cH_r,X_r,Y_r)\ge\Phi(S,X_r,Y_r)> e$, which is a contradiction, proving the lemma.
\end{proof}

\section{Making most of a random graph}
\label{sec:maker}

In this section we  prove Theorem~\ref{thm:makeGnp} by analysing the
following strategy of Maker.

\subsection*{Maker's strategy: degree and neighbourhood games}
Given $p$, $b$ and $\eps>0$, let $\delta=10^{-6}\eps^2$.  Maker
generates $\Gamma=G(n,p)$ flipping a $p$-biased coin for each edge
that she chooses to query as the game evolves.  Maker considers two
auxiliary games: the degree game and the neighbourhood game.  In odd
turns Maker plays the \emph{degree game}, which is a
\SB{1,2b,E(K_n),n,\ell_{\mathrm{D}},n} game with
$\ell_{\mathrm{D}}=\delta p n$. The underlying multihypergraph
$\cH^{\mathrm{D}}$ of the degree game is on vertex set~$E(K_n)$, with
hyperedges
$\partial_{K_n}(v):=\left\{uv\colon u\in[n]\setminus\{v\}\right\}$ for
each $v\in[n]$. The actions of the Ghost in this game are simply to
haunt any new member of $E(K_n)$ which Maker chooses to try to claim
in the real game and which turns out not to be in $\Gamma$ (the coin
comes up tails) or was claimed by Maker herself in the neighbourhood
game.  In particular, in this game the underlying
hypergraph~$\cH^{\mathrm{D}}$ does not change over time. The simulated
Breaker in this game simply claims the same subset of $E(K_n)$ in each
turn that the real Breaker claims in his two turns between successive
odd turns of Maker.

In even turns, Maker plays the \emph{neighbourhood game}, which we
will construct as a SpookyBalance game with parameters
$(1,2b,E(K_n),n,\ell_{\mathrm{N}},p^2n^2)$ with
$\ell_{\mathrm{N}}=\delta p^3n^2$. In this game, we start with an
underlying multihypergraph $\cHn_0$ on $E(K_n)$ consisting of $n$
empty hyperedges, one for each vertex in $[n]$. The actions of the Ghost
are as follows. First, in each round $r\ge 1$ of the neighbourhood
game (which corresponds to round $2r$ of the real game), the Ghost
creates $\cHn_r$. For this, it looks at each vertex $v\in[n]$ in the
real game, immediately before Maker's $2r$-th turn in the real game,
and adds to the hyperedge of $\cHn_{r-1}$ corresponding to~$v$ any
edge $xy$ of $K_n$ such that $xv$ and $yv$ have been claimed so far by
Maker in the real game, and such that~$xy$ is unclaimed in the real
game and is unhaunted in either game.  Second, the Ghost will haunt
any new member of $E(K_n)$ which Maker chooses to try and claim in the
real game and which turns out not to be in~$\Gamma$ or was claimed by
Maker herself in the degree game.  Again, the simulated Breaker in
this game simply claims the same subset of $E(K_n)$ in each turn that
the real Breaker claims in his two turns between successive even turns
of Maker.

In both games, Maker plays according to the potential function strategy, which
is a winning strategy for her by Lemma~\ref{lem:SB}. In both even and odd turns,
Maker claims the edge in the real game that she claims in the corresponding
auxiliary game.

We now prove that against any strategy of Breaker, when $\Gamma$ is drawn from the distribution $G(n,p)$, with high probability Maker wins the real game, proving Theorem~\ref{thm:makeGnp}. As mentioned earlier, we will analyse the point in the game when this strategy finishes: that is, all edges are either claimed or haunted. While Maker continues to play after this point (arbitrarily) the edges she claims later do not affect the conclusion of Theorem~\ref{thm:makeGnp}. We refer to this time (which is roughly round $p\binom{n}{2}$) as the time that the strategy finishes.

\begin{proof}[Proof of Theorem~\ref{thm:makeGnp}]
  Given $\eps\in (0,1)$, let $\delta=10^{-6}\eps^2$. Given
  $p\ge100\delta^{-1}n^{-1/2}$, suppose $b\le 10^{-6}\delta^3
  p^{-1}$. Fix a strategy of Breaker. We will reveal a graph $\Gamma$
  drawn from the distribution $G(n,p)$ as the game progresses: Maker's
  actions will reveal edges and non-edges of $\Gamma$. Breaker's
  actions will not reveal any edges or non-edges.  That is, when an
  edge~$xy\in E(K_n)$ is queried by Maker, we flip a $p$-biased coin
  to decide whether or not it is in~$\Gamma=G(n,p)$, but for the edges
  claimed by Breaker, we do not do this and postpone this decision to
  the end of the game.

 Before we analyse the result of the game, we pause to observe two properties $\Gamma$ satisfies with high probability. First, for each $v\in[n]$ we have $\deg_\Gamma(v)=(1\pm\delta)pn$. Second, for each $v\in [n]$ and $S\subseteq N_\Gamma(v)$ with $|S|\ge pn/2$, we have $e\big(G[S]\big)=(1\pm\delta)p|S|^2/2$. The first is a straightforward consequence of Lemma~\ref{lem:chernoff}, while to show the second, we begin by fixing $v$ and revealing $N_\Gamma(v)$. Assuming $\deg_\Gamma(v)\le(1+\delta)pn$, we next choose $S\subset N_\Gamma(v)$ with $|S|\ge pn/2$. We reveal the edges of $\Gamma$ in $S$, and by Lemma~\ref{lem:chernoff}, with input $\delta/2$, we find $(1\pm\delta)p|S|^2/2$ edges with probability at least $1-2\exp\big(-\tfrac{\delta^2p^3n^2}{100}\big)\ge 1-2\exp(-10pn)$, where the final inequality is by choice of $p$. Taking the union bound over the at most $2^{2pn}$ choices of $S$, we see that with probability at least $1-2\exp(-pn)$ all subsets $S$ of $N_\Gamma(v)$ with $|S|\ge pn/2$ satisfy $e\big(G[S]\big)=(1\pm\delta)p|S|^2/2$. Taking finally the union bound over choices of $v$, and since with high probability indeed $\deg_\Gamma(v)\le(1+\delta)pn$ for each $v$, we conclude that indeed the second property of $\Gamma$ also holds with high probability.
 
  Now let Maker play the above strategy. By Lemma~\ref{lem:SB}, she wins both auxiliary games: the 
  degree game and the neighbourhood game.   To see this, we simply need to observe that (respectively) the lower 
  bounds on $M$ and $\ell$ for given $m$, $b$ and $e$ are by the choice of
  $p$ satisfied. This is easy to check for~$M$, and we have
  \[\ell_{\mathrm{D}}=\delta p n\ge \tfrac{10b}{2b+1}\sqrt{\tfrac{n\log n}{2b+1}}\quad\text{and}\quad \ell_{\mathrm{N}}=\delta p^3n^2\ge \tfrac{10b}{2b+1}\sqrt{\tfrac{p^2n^2\log n}{2b+1}}\,,\]
  where in the latter we observe that at most $p^2n^2$ edges are in any
  vertex neighbourhood since no vertex neighbourhood has size more than
  $(1+\delta)pn$. Note that these inequalities hold with a great deal of
  room to spare; we only need $p\gg(\log n)^{1/3}n^{-2/3}$.

Our aim now is to show that Breaker claims only a  tiny fraction of edges
of $K_n$ at any given
vertex, and in any given vertex neighbourhood. We will see that this implies
that Maker claims most of the edges of the random graph $\Gamma$ in each set.
  
  For each $r\ge 1$, let $G^m_r$ be the graph of edges claimed by Maker
  after round $r$, and $G^b_r$ the graph of edges claimed by Breaker. It is
  convenient to view $(G^m_r,G^b_r)$ as a random process. We introduce one
  new quantity for the analysis. For $v\in[n]$, we let
  $S_0(v)=\emptyset$. For each $r$, if in round $r$, Maker claims the edge
  $uv$ and we have $\big|N_{G^m_{r-1}}(v)\cap N_{G^b_{r-1}}(u)\big|\ge
  \sqrt{\delta}pn$, then we set $S_r(v)=S_{r-1}(v)\cup\{u\}$, and otherwise
  we set $S_r(v)=S_{r-1}(v)$. In other words, we add $u$ to $S_r(v)$ when
  Maker claiming $uv$ adds exceptionally many Breaker-edges to the
  Maker-neighbourhood of $v$.
 We stop the process with failure at the first $r$ when for some $v$ we
 have either $\deg_{G^m_r}(v)\ge 2pn$ or $|S_r(v)|\ge
 8\sqrt{\delta}pn$. Otherwise we terminate the process when every edge of $K_n$ is
 claimed by either Maker or Breaker, or is haunted. That is, this is the point where Maker has queried containment in $\Gamma$ of every edge of $K_n$ not claimed by Breaker.  We
 now show that neither failure event is likely to occur.

 For each $r$ and $v$, let $T_r(v)$ be the set of vertices $u$ such that
 $\big|N_{G^m_{r-1}}(v)\cap N_{G^b_{r-1}}(u)\big|\ge \sqrt{\delta}pn$. In
 other words, $T_r(v)$ is the set of vertices which send exceptionally many
 Breaker-edges to the Maker-neighbourhood of $v$. Observe that if we have
 not stopped the process, we have $\big|N_{G^m_{r-1}}(v)\big|<2pn$ for each
 $v\in[n]$. Because Maker wins the degree game, it follows that
 $\big|N_{G^b_{r-1}}(u)\big|<4bpn+(2b+1)\delta p n<\delta n$. In
 particular, at most $2\delta pn^2$ edges in $G^b_{r-1}$ are incident to
 $N_{G^m_{r-1}}(v)$, so the Hand-Shaking Lemma gives $|T_r(v)|\le 4\sqrt{\delta}n$. Furthermore, observe
 that $T_{r-1}(v)\subseteq T_{r}(v)$ for each $r,v$. Observe that both
 $\deg_{G^m_r}(v)$ and $|S_r(v)|$ are given by a sum of independent
 Bernoulli random variables with probability $p$. In the former case,
 trivially at most $n-1$ random variables are summed, while in the latter
 case at most $|T_r(v)|\le 4\sqrt{\delta}n$ random variables are
 summed. Since these sums are thus stochastically dominated by binomial
 random variables with expectations respectively $pn$ and
 $4\sqrt{\delta}pn$, we conclude that by Lemma~\ref{lem:chernoff}, with
 high probability the process does not stop with failure. From now on we
 assume that the process does not stop with failure.

 Let $G^m$ and $G^b$ be the graphs of edges claimed by respectively Maker
 and Breaker when the strategy finishes, and let $S(v)$ be the final set
 $S_r(v)$. We will now show that $G^b$ contains only a tiny fraction of
 edges of $K_n$ at each vertex and in each $G^m$-neighbourhood. Observe that none of the edges of $G^b$ has been revealed when the process terminates. We already showed $\deg_{G^b}(v)<\delta n$ for each $v$, so that what remains is to show Breaker claims few edges in each $G^m$-neighbourhood.
 
 Fix a vertex $v$. We split the edges Breaker claims in $N_{G^m}(v)$ into three classes, according to when Breaker claims them. For a given pair $xy$ contained in the neighbourhood  $N_{G^m}(v)$ which is in $G^b$, let $r_{xy}$ denote the round of the real game when Maker claimed the later of $vx$ and $vy$. If $xy$ remains neither claimed nor haunted immediately before the next turn of Maker in the neighbourhood game, which is in round $2\big(\lfloor\tfrac{r_{xy}}{2}\rfloor+1\big)$ of the real game, then $xy$ is added to the hyperedge corresponding to $v$ in the neighbourhood game. Thus, we let the first class of Breaker-edges be $xy$ which Breaker claimed in round $2\big(\lfloor\tfrac{r_{xy}}{2}\rfloor+1\big)$ or later of the real game; these are exactly the edges that the simulated Breaker claims in the hyperedge corresponding to $v$ in the neighbourhood game.
 
 The second class of edges are $xy$ which Breaker claimed in round $r_{xy}-1$ or before of the real game; these are edges which Breaker is `lucky' to claim, because he does not know they will eventually be in the Maker-neighbourhood of $v$.
 
 Finally, the third class of edges are the rest: edges which Breaker claims in round $r_{xy}$ and (if $r_{xy}$ is even) $r_{xy}+1$ of the real game, after he knows that $xy$ is in the Maker-neighbourhood of $v$ but before Maker's next turn in the neighbourhood game.
 
 To deal with the first class of edges, we recall that Maker wins the neighbourhood game, and that Maker does this by claiming at most all the edges of $\Gamma$ in $N_{G^m}(v)$. Since $\big|N_{G^m}(v)\big|\le \big|N_\Gamma(v)|\le 2pn$, the number of edges Maker claimed in $N_{G^m}(v)$ is $M\le 2p(2pn)^2/2=4p^3n^2$. Letting $B$ denote the number of first-class edges, we have, by definition of Maker winning the neighbourhood game, $\tfrac{M+B}{1+2b}-\ell_{\mathrm{N}}\le M$. Rearranging, we get
 \[B\le 2bM+(1+2b)\ell_{\mathrm{N}}\le 8bp^3n^2+3b\delta p^3n^2\le 10bp^3n^2=10^{-5}\delta^3p^2n^2\,.\]
 
 We split the second class of edges into two parts: those using a vertex of $S(v)$ and those not (which must still use a vertex of $N_{G^m}(v)$, which has size at most $2pn$). A vertex of $S(v)$ is in at most $|N_{G^m}(v)|\le 2pn$ second-class edges, while a vertex which is not in $S(v)$ is by definition in at most $\sqrt{\delta}pn$ second-class edges. Summing, and using our bound on $|S(v)|$, the total number of second-class edges is at most $8\sqrt{\delta}pn\cdot 2pn+2pn\cdot\sqrt{\delta}pn=18\sqrt{\delta}p^2n^2$.
 
 Finally, for the third class of edges, we note that the possible values of $r_{xy}$ are the rounds at which Maker claimed an edge at $v$, of which there are at most $2pn$, so the number of third-class edges is at most all of the $4bpn\le 10^{-5}\delta^3 n$ edges that Breaker claimed in the at most $4pn$ rounds of the form $r_{xy}$ and $r_{xy}+1$.
 
 Putting these bounds together, we conclude the total number of Breaker-edges in $N_{G^m}(v)$ is at most
 \[10^{-5}\delta^3p^2n^2+18\sqrt{\delta}p^2n^2+10^{-5}\delta^3 n\le 20\sqrt{\delta}p^2n^2\,,\]
 where we used $p^2n\le 1$ for the inequality.
%
%
 
 Finally, we reveal the edges of $\Gamma$ in $G^b$: recall that these were not revealed previously. The point here is that all edges of $\Gamma$ not claimed by Breaker must have been claimed by Maker. Using
 Lemma~\ref{lem:chernoff} and the union bound, we see that with high
 probability at most $2\delta pn$ edges of $\Gamma$ appear in $G^b$ at
 any vertex $v\in[n]$, and hence
 $\deg_{G^m}(v)\geq(1-3\delta)pn\geq(1-\eps)pn$.  Again using
 Lemma~\ref{lem:chernoff} and the union bound, with high probability
 Breaker claimed at most $40\sqrt{\delta}p^3n^2$ edges of $\Gamma$ in
 $N_{G^m}(v)$. By the second property of $\Gamma$,
 the total number of edges of $\Gamma$ in $N_{G^m}(v)$
 is at least $(1-\delta)p(1-3\delta)^2p^2n^2/2$, of which at most
 $40\sqrt{\delta}p^3n^2$ are claimed by Breaker. The remaining at least
 $(1-\eps)p^3n^2/2$ edges are claimed by Maker, giving the desired
 edges claimed by Maker in $N_\Gamma(v)$. 
\end{proof}

Let us finish this section with a few remarks about simple modifications to
the setup considered in Theorem~\ref{thm:makeGnp}.  First, it is easy to
check that this strategy also works if $K_n$ is replaced by a typical
random graph $G(n,q)$ for some $q>p$.  In that case the random graph
$G(n,p)$ is obtained from $G(n,q)$ by keeping edges independently with
probability $p/q$ and we require Breaker's bias to satisfy
$b\ll q/p$.  Secondly, it is similarly easy to check that this proof goes
through if Breaker, rather than Maker, plays first. To see this, observe
that removing any single edge from Maker's graph (corresponding to Maker
losing her first turn) only affects the number of edges Maker claims at a
vertex or in a vertex neighbourhood by one; this off-by-one error is
absorbed by the slack in the constant choices.

\section{Universal graphs}
\label{sec:resilience}

The following resilience result for random graphs together with
Theorem~\ref{thm:makeGnp} immediately implies Corollary~\ref{cor:makeuniv}.

\begin{theorem}\label{thm:universalsubgraph}
For all $D,\,\Delta\ge 2$ there exist $\eps>0$ and $C$ such that if
  $p\ge C\big(\tfrac{\log n}{n}\big)^{1/D}$, then the random graph $\Gamma=G_{n,p}$ asymptotically almost
  surely has the following property.
  
  Let $G$ be any spanning subgraph of $\Gamma$ such that $\deg_G(v)\ge(1-\eps)pn$ and $e_G\big(N_{\Gamma}(v)\big)\ge(1-\eps)p^3n^2/2$. Then $G$ contains all graphs $H$ on at most $n$ vertices with $\Delta(H)\le D$, and all graphs $H$ on at most $n$ vertices with $\Delta(H)\le\Delta$ and degeneracy at most $\tfrac{D-1}{2}$.
\end{theorem}
\begin{proof}
 Given $D$ and $\Delta$, without loss of generality we assume $\Delta\ge D$ and we set $r=8\Delta$. We set $d=\tfrac12$. Let $\eps'>0$ and $C$ be returned by Lemma~\ref{lem:simpBL}. We set $\eps=10^{-6}r^{-2}(\eps')^3$. Let $p\ge C\big(\tfrac{\log n}{n}\big)^{1/D}$.
 
 Partition $[n]$ into $r$ parts $V_1,\dots,V_{r}$, in size order,
 differing in size by at most one. We now generate
 $\Gamma=G(n,p)$. Using Lemma~\ref{lem:chernoff}, we can easily check
 that with high probability, the following all hold. First, for each
 $v\in V(\Gamma)$ we have $\deg(v)=(1\pm\eps)pn$ and $N(v)$ contains
 $(1\pm\eps)p^3n^2/2$ edges. Second, for each $i\in[r]$ and $v\not\in
 V_i$ we have $\deg(v;V_i)=(1\pm\eps)p|V_i|$. Third, for each $i\neq
 j\in[r]$, and each $v\not\in V_i$, the pairs $(V_i,V_j)$, and $\big(N(v;V_i),V_j\big)$ are $(\eps,1,p)$-regular, and if $v\not\in V_j$ in addition $\big(N(v;V_i),N(v;V_j)\big)$ is $(\eps,1,p)$-regular. Suppose now that $\Gamma$ satisfies this good event and in addition the good event of Lemma~\ref{lem:simpBL} with the inputs as above.
 
 Let $G\subset\Gamma$ be such that $\deg_G(v)\ge(1-\eps)pn$ and
 $e_G\big(N_{\Gamma}(v)\big)\ge(1-\eps)p^3n^2/2$. Then $G$ and
 $\Gamma$ differ by at most $2\eps pn$ edges at each vertex, and at
 most $\eps p^3n^2$ edges in each vertex neighbourhood. In particular,
 for any $i$ and $S\subset V_i$, at most $2\eps pn|S|$ edges were
 removed from the vertices of $S$, and so that for any $j\neq i$ and
 $T\subset V_j$ with $|T|\ge\eps'|V_j|$, the pair $(S,T)$ has density
 at least $p/2$. It follows that $(V_i,V_j)$ and
 $\big(N_\Gamma(v;V_i),V_j)$ are $(\eps',\tfrac12,p)$-regular in $G$
 for each $i\neq j$ and $v\not\in V_i$. Finally, since
 $\big(N_\Gamma(v;V_i),N_\Gamma(v;V_j)\big)$ is $(\eps,1,p)$-regular
 in $\Gamma$ for any $i\neq j$ and $v\not\in V_i\cup V_j$, and since
 at most $\eps p^3n^2$ edges were removed from $N_\Gamma(v)$ to form
 $G$, also the pair $\big(N_\Gamma(v;V_i),N_\Gamma(v;V_j)\big)$ is
 $(\eps,\tfrac12,p)$-regular in $G$. We conclude that $G$, with the
 partition $V_1,\dots,V_r$, meets the conditions of
 Lemma~\ref{lem:simpBL}.
 
 If $H$ satisfies $v(H)\le n$ and $\Delta(H)\le D$, we can add
 isolated vertices if necessary to obtain a graph with exactly $n$
 vertices. By Lemma~\ref{lem:nicepartition} with input $Y=\emptyset$
 we can let $X_1,\dots,X_r$ be a partition of $V(H)$ into independent
 sets, in size order, differing in size by at most one, and hence
 $|V_i|=|X_i|$ for each $i$. Now Lemma~\ref{lem:simpBL} states that
 $H\subset G$, as desired.
 
 Finally, if $H$ satisfies $v(H)\le n$, and $\Delta(H)\le \Delta$, and $H$ has degeneracy at most $\tfrac{D-1}{2}$, then again we can add isolated vertices to $H$ if necessary to obtain an $n$-vertex graph. Let $Y$ be the set of vertices in $H$ with degree at most $2D$. Since $H$ has less than $Dn$ edges, we have $|Y|\ge\tfrac{n}{2D+1}$. Now let $X_1,\dots,X_r$ be a partition of $V(H)$ into independent sets, in size order, obtained by applying Lemma~\ref{lem:nicepartition} with input $Y$. We have $|V_i|=|X_i|$ for each $i$, and $|X_i\cap Y|\ge\tfrac{1}{2(2D+1)}|X_i|$ for each $i$. Thus again Lemma~\ref{lem:simpBL} states that $H\subset G$, as desired.
\end{proof}

\section*{Acknowledgement}

We are grateful for the helpful suggestions of two anonymous referees.

Data access statement: This article does not report any underlying research data.

\bibliographystyle{amsplain}
\bibliography{extracted-mrefed, additional}
\end{document}